\documentclass[epsf]{article}
\usepackage{amsfonts}
\usepackage{subfigure} 
\usepackage{epsfig}
\usepackage{amsmath, amsthm}
\newcommand{\be}{\begin{eqnarray}}     	\newcommand{\ee}{\end{eqnarray}}

\newcommand{\vol}{\mathrm{Vol}}
\newcommand{\diam}{\mathrm{diam}}
\newcommand{\sing}{\mathrm{Sing}}
\newcommand{\ric}{\mathrm{Ric}}
\newcommand{\supp}{\mathrm{supp}}
 
\title{Convergence of Kahler-Einstein orbifolds}
\author{Natasa Sesum}
\theoremstyle{plain}
\newtheorem{dummy}{Dummy}

\theoremstyle{definition}

\newtheorem{lemma}[dummy]{Lemma}
\newtheorem{theorem}[dummy]{Theorem}
\newtheorem{proposition}[dummy]{Proposition}
	
\newtheorem{definition}[dummy]{Definition}

\newtheorem{claim}[dummy]{Claim}

\begin{document}

\maketitle

\begin{abstract}

We proved the convergence of a sequence of 2 dimensional comapct
Kahler-Einstein orbifolds with rational quotient singularities and
with some uniform bounds on the volumes and on the Euler
characteristics of our orbifolds to a Kahler-Einstein 2-dimensional
orbifold. Our limit orbifold can have worse singularities than the
orbifolds in our sequence. We will also derive some estimates on the
norms of the sections of plurianticanonical bundles of our orbifolds
in the sequence that we are considering and our limit orbifold.

\end{abstract}

\begin{section}{Introduction}

The main goal of this paper is to try to generalize some theorems in
\cite{tian1990} that hold for the sequences of smooth Kahler-Einstein
complex surfaces to the sequences of $2$ dimensional Kahler-Einstein
orbifolds.

\begin{definition}
A complex orbifold is a complex manifold $M$ of dimension $n$ whose
singularities are locally isomorphic to quotient singularities $C^n/G$
for finite subgroups $G\subset GL(n,C)$. 
We say that $g$ is a Kahler metric on a complex orbifold $M$ if $g$ is
Kahler in the usual sense on the nonsingular part of $M$ and whenever
$M$ is locally isomorphc to $C^n/G$, we can identify $g$ with the
quotient of a $G$ - invariant Kahler metric defined near $0$ in $C^n$.
\end{definition}

Let $\mathcal{J}_n$ be the collection of all complex surfaces of the
form $CP^2 \# n\overline{CP^2}$ (i.e. all complex structures on
$CP^2\#\overline{CP^2}$).

The following theorem was proved by Tian in \cite{tian1990}.

\begin{theorem}[Tian]
\label{theorem-theorem_tian}
Let $\{(M_i,g_i)\}$ be the sequence of compact Kahler-Einstein
manifolds , where $(M_i,g_i)$ is a comact complex surface in
$\mathcal{J}_n$ ($5 \le n \le 8)$. Then by taking a subsequence if
necessary, we may assume that $(M_i,g_i)$ converge to a
Kahler-Einstein manifold $(M_{\infty}\backslash \sing(M_{\infty}),
g_{\infty})$, where $(M_{\infty}, g_{\infty})$ is a connected
Kahler-Einstein orbifold and $\sing(M_{\infty})$ is the finite set of
singular points of $M_{\infty}$.
\end{theorem}

In this paper $\{(M_i, g_i)\}$ will always denote a sequence of $2$
dimensional Kahler-Einstein orbifolds with positive first Chern class
and with rational singularities. Let $C_1$, $C_2$ be arbitrary
constants. Define $\mathcal{A}(C_1,C_2)$ to be the following set

$\mathcal{A}(C_1,C_2) = \{(M,g)$ orbifold as above $\:\:|\:\: C_1 \le
\vol_g M,\:\: |\chi(M)| \le C_2\}$.

\begin{definition}
\label{definition-definition_convergence}
We will say that a sequence of Kahler-Einstein orbifolds $\{M_i,g_i\}$
converges to a Kahler-Einstein orbifold $(M_{\infty}, g_{\infty})$ if
the sequence converges in Hausdorff topology to a conneceted Einstein
orbifold $M_{\infty}$ and if:

Let $\{p_i\}_{1\le i \le N}$ be singular points of $M_{\infty}$. If
$G_{\infty} = M_{\infty}\backslash \bigcup \{p_i\}$, then $G_{\infty}$
has a $C^{\infty}$ Einstein metrics $g_{\infty}$ and there are
$C^{\infty}$ embeddings $F_i : G_{\infty} \to M_i$ for $i$
sufficiently large, such that on every compact set $K$ of
$G_{\infty}$:

\begin{enumerate}

\item
$F_i^*g_i$ converge to $g_{\infty}$ uniformly on $K$.

\item 
$(F_i^{-1})^*\circ J_i\circ F_i*$ converge to $J_{\infty}$
uniformly on $K$, where $J_i$, $J_{\infty}$ are the almost complex
structures of $M_i$, $M_{\infty}$, respectively.

\end{enumerate}

Moreover, each singular point $p_i$ has a neighbourhood which is
homeomorphic to a cone on a spherical space form
$C(S^{n-1}/\Gamma)$. If the metric $g_{\infty}$ is lifted to
$B^n\backslash{0}$ via $\Gamma$, then there is a $\Gamma$-equivariant
diffeomorphism $\phi: B^n\backslash {0} \longrightarrow B^n\backslash
{0}$ such that $\phi^*g_{\infty}$ extends smoothly over $0$ to a
smooth Einstein metric on $B^n$.\\ \\ 

\end{definition}

The main purpose of the paper is to prove the following theorem

\begin{theorem}
\label{theorem-theorem_main_theorem}
Let $\{(M_i, g_i)\}$ be the sequence of compact Kahler-Einstein $2$
dimensional orbifolds in $\mathcal{A}(C_1,C_2,C_3)$. Then by taking a
subsequence if necessary, we may assume that $(M_i, g_i)$ converge to
a Kahler-Einstein orbifold $(M_{\infty}, g_{\infty})$ with a finite set
of singular points.
\end{theorem}

The organization of this paper is as follows. In section $2$ we will
give some preliminaries and we will prove our main theorem about the
convergence of a sequence of Kahler-Einstein $2$ dimensional
orbifolds. The argument is based on the arguments of Tian in
\cite{tian1990} and the arguments of Bando, Kasue and Nakajima in
\cite{bando1989}. In section $3$ we will give some applications of our
main theorem to the convergence of a sequence of global holomorphic
sections of plurianticanonical bundles $H^0(M_i, K_{M_i}^{-m})$ to a
global holomorphic section of $K_{M_{\infty}}^{-m}$.

I would like to thank my advisor Gang Tian for his support, guidance
and bringing this problem to my attention. I would like also to thank
Jeff Viaclovsky for helpful suggestions and discussions.

\end{section}

\section{Compactness argument for the sequence of orbifolds}

Let $(M_{i},g_{i})$ be the sequence of 2-dimensional Kahler-Einstein
orbifolds,with positive first Chern class and with rational
singularities (meaning that for each orbifold group $G < SU(2)$), such
that $\diam(M_{i})\leq C_1,\:\: C_2 \leq \vol_{g_i}, | \chi(M_i) |
\leq C_3,\:\: \forall i$ for some uniform constants $C_1,C_2,C_3$. We
may assume that $\ric(g_i)=w_{g_i}$.

\begin{theorem}
\label{theorem-order-singularity}
The orders of singularities of all $M_{i}$ are uniformly bounded,
i.e. $\exists C$ such that $| \gamma_{i}^{j}| \leq C,\:\:\: \forall i$
and $\forall j \in \{1,\dots\, N_i\}$, where $\gamma_{i}^{j}$ is the
orbifold group associated to the singularity $p_{i}^{j} \in M_i$ and
$N_i$ is the number of isolated singularities of $M_{i}$.
\end{theorem}

\begin{proof}
Take $p_{i}^{j}\in M_{i}$, and let $r \ge 0$ be such that
$B(p_{i}^{j})\cong \Delta ^{r}/\gamma_{i}^{j}$, where $\Delta^{r}$ is
a disk in $C^2$. Then:

\begin{equation}
\label{equation-equation_order_singularity1}
\vol_{g_i}(B(p_{i}^{j},r)) =
\frac{\vol_{\tilde{g_i}}(\Delta^{r})}{|\gamma_{i}^{j}|}
\end{equation}

From the lower bound for the Sobolev constant for orbifolds (see
\cite{nakagawa1993}) which is uniform in $i$ in the case of our
sequence of orbifolds we have:

\begin{equation}
\label{equation-equation_order_singularity2}
\vol_{g_i}B(p_{i}^{j}) \geq Cr^4
\end{equation}

for some uniform constant C and small values of $r$. Upstairs on
$\Delta^r$ we have nonnegative Ricci curvature so Bishop-Gromov
comparison principle implies that 

\begin{equation}
\label{equation-equation_euclidean}
\vol (\Delta^r) \le w_n r^4
\end{equation}

where $w_n$ is the euclidean constant. Equations
\ref{equation-equation_order_singularity1} and
\ref{equation-equation_euclidean} imply downstairs on $M_i$:

$$\vol_{g_i}B(p_i^j,r) \le \frac{Cr^4}{|\Gamma_i^j|}$$

Combining this with equation
\ref{equation-equation_order_singularity2} for small $r$ gives us 

$$|\gamma_{i}^{j}| \leq C$$

for some uniform constant C, $\forall i,j$.

\end{proof}

\begin{theorem}
\label{theorem-theorem_number_singularities}
There exists a uniform bound on the number of singularities,
i.e. $\exists C$, s.t. $N_i \leq C, \forall i$.
\end{theorem} 

\begin{proof}

\begin {equation}
\label{equation-equation_integral}
\int_{M}c_{1}^2(M) - 2c_2(M) = 3\sigma(\tilde{M}) + 2\sum_{p\in
\sing(M)}(e(E_{p}) - \frac{1}{G_{p}})
\end{equation}

where $E_{p}$ is the exceptional divisor of the minimal
desingularization $\pi: \tilde{M}\to M$ and $E_p =
\bigcup_{j=1}^{k_p}C_{pj}$, where p is a singularity, and $k_p$ is the
order of singularity at $p$. Since $p$ is a rational singularity,
$C_{pj}$ are rational curves. From algebraic geometry it follows that
we can choose $\pi$ to be a composition of consequtive blow-ups, such
that $C_{pj}^2 \leq -2 \:\:\: \forall p,\:\: \forall j$ and
$C_{pj}C_{pk} = 1 \:\:\: \forall j \neq k$ and better, s.t. no 3
distinct $C_{pj}$ meet. It now follows that the number of pairs
$C_{pj}, C_{pk}$ such that $C_{pj}C_{pk}=1$ is less or equal than $k_p
- 1$ where $k_p \leq C$ and $C$ is taken from theorem
\ref{theorem-order-singularity}.  Topologically, $E_p$ is a connected
sum of $k=k_p$ copies of $CP^1$ and therefore $\chi(E_p) = 2$.

\begin{equation}
\sigma(\tilde{M})= \frac{1}{3}(c_1(\tilde{M})^2 - 2c_2(\tilde{M}) \nonumber
\end{equation}

\begin{equation}
c_1(\tilde{M}) = \pi^{*}c_1(M) - E_1 - \dots - E_N \nonumber
\end{equation}

\begin{equation}
c_1(\tilde{M})^2 = \pi^{*}c_1(M)^2 + \sum_{j=1}^{j=N} E_j^2 \nonumber
\end{equation}

\begin{equation}
E_i^2 = \sum_{j}C_{ij}^2 + \sum_{j \neq k} C_{ij}C_{ik} \leq -2k_i +
2(k_i -1) = -2 \nonumber
\end{equation}

\begin{equation}
\label{equation-c1}
c_1(\tilde{M})^2 \leq \pi^*c_1(M)^2 - 2N
\end{equation}

From equation \ref{equation-equation_integral} we get that:

$$ \int_M c_1^2(M) - 2c_2(M) \leq 3\int_{\tilde{M}}(c_1^2(\tilde{M}) -
2c_2(\tilde{M})) + 4N - \frac {2}{C}N$$

where C is taken from theorem \ref{theorem-order-singularity}. From
equation \ref{equation-c1}:

$$2N\leq 2(c_2[\tilde{M}] - c_2[M])\leq \pi^*c_1[M]^2 - c_1[M]^2 + 2N
- \frac{2N}{C}$$

From above we get (since $c_1^2[M] \le \bar{C}$ is uniformly
bounded):

$$N \le \frac{1}{2}C\cdot 2C_1 = \tilde{C}$$

for some uniform constant $\tilde{C}$. 
 
\end{proof}  

In the proof we wrote $E_i$ instead of $E_{p_i}$, where $p_i$ is a
singular point ( similarly $C_{ij}$, $k_i$ and $N_i$ are related to a
point $p_i$).

The following theorem can be found in \cite{anderson1989}.

\begin{theorem}[Anderson]
\label{theorem-theorem_anderson}
There is a constant $C = C(n, c_S)$ and $\epsilon_0 =
\epsilon_0(n,c_S)$ such that if $B(t)$ is a geodesic ball of radius
$t$ in $M$ (where $(M,g)$ is a Kahler-Einstein surface with $\ric(g) =
w_g$ and

$$\int_{B(t)}|R|^{\frac{n}{2}} dV < \epsilon_0$$

then

$$\sup_{B(\frac{t}{2})}|R| \le C\cdot
\frac{1}{t^2}(\int_{B(t)}|R|^{\frac{n}{2}})^{\frac{2}{n}}$$
\\
\end{theorem}

Anderson proved theorem \ref{theorem-theorem_anderson} for smooth
manifolds but it holds for our orbifolds with isolated singularities
as well. Let $(M,g)$ be a Kahler-Einstein orbifold with rational
isolated singularities. The following inequalities for the Laplacian
of the curvaure tensor hold for $M$ in a weak sense (in the sense of
distribution):

$$\Delta R = R * R + R * \ric + P^2(\ric)$$

where $A * B$ denotes a linear combination of tensors $A$, $B$
obtained by contracting $A$, $B$ with the metric $g$ and $P^2(\ric)$
is a linear combination of second covariant derivatives of the Ricci
tensor. In particular, one obtains

$$|\Delta R| \le c_1|D^2\ric| + c_2|R|^2$$

where $c_1$ and $c_2$ are constants depending on dimension.  Furthermore

$$\langle \Delta R,R\rangle + |DR|^2 = \frac{1}{2}\Delta|R|^2 =
|R|\Delta|R| + |d|R||^2$$

We have that $|d|R||^2 \le |DR|^2$. By Schwartz inequality applied to
$\langle \Delta R,R\rangle$ we get

$$\Delta|R| + c_1|D^2\ric| + c_2|R|^2 \ge 0$$

Since $\ric(g) = w_g$ we have $D\ric = 0$ and therefore 

$$\Delta|R| + c_2|R|^2 \ge 0$$

holds for orbifold $M$ in a weak sense.

In his paper \cite{nakagawa1993} Nakagawa proved that Sobolev
inequality holds on orbifolds with the lower bound on Ricci curvature,
i.e.

$$||f||_4 \le \frac{1}{c_S}||\nabla f||_2 + \vol^{-1}||f||_2$$

for any Lipshitz function for $M$. $c_S$ is a Sobolev constant that is
for our sequence of orbifolds uniform (because of uniform bounds on
our sequence of orbifolds $\{(M_i,g_i)\}$ specified at the beginning
of the paper).
                                                         
Sobolev inequality holds for orbifolds and therefore the Moser
iteration argument (as in \cite{anderson1989}) gives us that $\exists$
$C = C(n,C_s)$, $\epsilon_0 = \epsilon_0(n,C_s)$, where $C_s$ is a
Sobolev constant, such that if B(t) is a geodesic ball of radius t in
M (M is an orbifold) and if $\int_{B(t)}|R|^{\frac{n}{2}}dV \le
\epsilon_0$, then:
$$\sup_{B(\frac{t}{2})}|R| \leq
\frac{C}{t^2}(\int_{B(t)}|R|^{\frac{n}{2}})^{\frac{2}{n}}$$

Since there is a uniform bound on the number of singularities, by
taking a subsequence of orbifolds, we may assume that each $M_i$ has
$S$ singular points $\{p_j^i\}_{1\le j \le S}$.  Now following the
arguments in \cite{tian1990} for a sequence of smooth surfaces, we can
conclude that in the case of a sequence of orbifolds there exists a
subsequence $(M_i,g_i)$, such that $M_i \backslash \{ \{
x_{i_\beta}\}_{1 \leq \beta \leq N} \cup \{p_{j}^{i}\}_{1 \leq j \leq
S}\}$ converge to a Kahler Einstein manifold $(M_{\infty},g_{\infty})$
in the sense of a definition
\ref{definition-definition_convergence}. Since a distance functions on
$M_i \times M_i$ converge to a Lipshitz function $\rho_{\infty}$, the
same argument as in \cite{tian1990} shows that we can attach finitely
many points $x_{\infty 1} \dots x_{\infty N}$ and $p_{\infty 1}, \dots
, p_{\infty S}$ to $M_{\infty}$ such that we get a complete metric
space. From \cite{tian1990} we know that $x_{\infty\beta}$ for $1 \leq
\beta \leq N$ are the orbifold points (we get this points in the limit
process as a result of concentrating a curvature of $(M_i,g_i)$ at
smooth points $\{x_{i\beta}\}_{1\le\beta\le N}$).

To finish the proof of therem \ref{theorem-theorem_main_theorem} we
only need to check that $\{p_{i\infty}\}_{1\le i\le S}$ are the
orbifold points of $M_{\infty}$. These points come from singular
points $p_j^i$ of our orbifolds $M_i$.

Let $p_{\infty i} = p$ and look at $B = B(p,t_0)$, a ball in a
complete metric space $M_{\infty}$. We want to show that the ball $B$
satisfies the theorem proved in \cite{bando1989} (we will state it
below), since then we will be able to conclude that $p$ is an orbifold
point. Without loss of generality assume that $p$ is the only singular
point of $M_{\infty}$.

\begin{theorem}[Bando, Kasue, Nakajima]
\label{theorem-theorem_Bando}
Let $B = B(p,t_0)$ be ball in a complete, locally compact metric space
of length $\rho_{\infty}$. Suppose that $B\backslash\{p\}$ is locally
connected, i.e. for every open set U containing p, there exists an
open set V, containing p, s.t. $V \subset U$ and $V\backslash \{p\}$
is connected $C^{\infty}$ manifold with Einstein metric g, satisfying:

\begin{enumerate}

\item
\label{cond1}
$$\int_{B\backslash\{p\}}|R|^2 \le \infty$$

\item
\label{cond2}

$$(\int_{B\backslash\{p\}}|v|^4)^{\frac{1}{2}} \leq
S\int_{B\backslash\{p\}}|Dv|^2 \ \ \forall v \in C_{0}^{1}(B\backslash
\{p\})$$

\item
\label{cond3}

$$\vol B(p,t) \leq Vt^4 \ \ \   \forall \ \ 0 \le t \le t_0$$

\end{enumerate}

Then the metric $g_{\infty}$ extends smoothly to $B$ as an orbifold
metric.

\end{theorem}

Gauss-Bonnet theorem applied to our orbifolds $(M_i,g_i)$ gives us:

\begin{equation}
\label{equation-equation_gauss_bonnet}
\chi(M_i) + \sum_{i=0}^{N_i}\frac{1}{|\Gamma_i|} =
\frac{1}{8\pi^2}\int_{M_i}(|R_i|^2 - 4|\ric_i|^2 + \tau_i^2)
\end{equation}

where $\tau_i$ is a scalar curvature. Since our orbifolds are
Kahler-Einstein, the curvature integral in formula
\ref{equation-equation_gauss_bonnet} becomes $\int_{M_i}|R_i|^2$.
Euler characteristic bounds for $(M_i,g_i)$ and uniform bounds on the
number of singularities of $M_i$ give us uniform upper bound of LHS of
formula \ref{equation-equation_gauss_bonnet}. Therefore there exists
some constant $C$ (independent of $i$) such that

$$\int_{M_i}|R_i|^2 dV_{g_i} \le C$$

for some uniform constant $C$. In particular, it implies that $L^2$
integral of $||R_i||_{g_i}$ ic uniformly bounded from above by a
uniform constant.

By Fatou's lemma we now get that $\int_B |R|^2_{g\infty} \le C <
\infty$, i.e. the condition \ref{cond1} of theorem
\ref{theorem-theorem_Bando} is satisfied.

\begin{lemma}
\label{lemma-cond2}
Condition \ref{cond2} is satisfied for $B$.
\end{lemma}

\begin{proof}

Let $v \in C_{0}^{1}(B\backslash \{p\})$ and let $\supp(v) = K
\subset B\backslash\{p\}$. By the definition of convergence , there
exist diffeomorphisms $\phi_i$ from the open subsets of $M_i\backslash
\{p_i\}$ to the open subsets of $M\backslash\{p\}$ that contain $K$,
such that every diffeomorphism $\phi_i$ maps some compact subset
$K_{i}$ onto $K$, where $K_{i}$ is contained in $B(p_i, t_0)$, for
some sufficiently large $i$ (because of the uniform convergence of
metrics on compact subsets). We have that $\tilde{g_i}
=(\phi_{i}^{-1})^*g_i$ converge uniformly and smoothly on $K$ to
$g_{\infty}$.

Since $M_i$ is an orbifold, the Sobolev inequality holds, with a
constant that does not depend on $i$ (because of our uniform bounds on
the sequence $\{(M_i,g_i)\}$ as in theorem
\ref{theorem-theorem_main_theorem}). Let $F_i = \phi_{i}^* (v)$. Then,
$\supp F_i \subset K_i \subset B(p_i,t_0)\backslash \{p_i\}$. Let
$\{\eta_{i}^{k}\}$ be the sequence of cut-off functions, such that
$\eta_{i}^{k} \in C_{0}^{1}(B_i \backslash \{p_i\})$ and $\eta_{i}^{k}
\to 1 (k \to \infty) \ \ \forall i$, and:
$$\int_{B_i}|D\eta_{i}^{k}|^2 \to 0 \ \ (k\to\infty)$$

$\eta_{i}^{k} F_i$ is a function of compact support in
$B(p_i,t_0)$. Then by Sobolev inequality:

$$(\int_{B(p_i,t_0)}|\eta_{i}^{k} F_i|^4 dV_{g_i})^{\frac{1}{2}} \leq
C\int_{B(p_i,t_0)}|D(\eta_{i}^{k}F_i)|^2 dV_{g_i}$$

We can bound $F_i$ with some constant $C_i$, and therefore:

$$\int_{B(p_i,t_0)}|D(\eta_{i}^{k}F_i)|^2 \leq
C(\int_{(B(p_i,t_0))}|D\eta_{i}^{k}|^2 C_i +
\int_{(B(p_i,t_0)}|DF_{i}|^2(\eta_{i}^{k})^2)$$

Let $k$ tend to $\infty$. Then we get:

$$(\int_{B(p_i,t_0)\backslash \{p_i\}}|F_{i}|^4)^{\frac{1}{2}} \leq
C\int_{B(p_i,t_0)\backslash \{p_i\}}|DF_i|^2$$

Since $\supp F_i \subset K_i$, after changing the coordiantes, via map
$\phi_i$ we get:

$$(\int_{K}|v|^4 dV_{\tilde{g_i}})^{\frac{1}{2}} \leq C \int_{K}|Dv|^2
dV_{\tilde{g_i}}$$

$\tilde{g_i}$ converges uniformly on K to $g_{\infty}$, so letting $i$
tend to $\infty$ in the above inequality, keeping in mind that $\supp v
= K$, we get that:

$$(\int_{B\backslash\{p\}}|v|^4 dV_{g_{\infty}})^{\frac{1}{2}} \leq
C\int_{B\backslash\{p\}}|Dv|^2 dV_{g_{\infty}}$$

\end{proof}	

\begin{lemma}
\label{lemma-cond3}
$\vol(B(p,t)) \leq Ct^4$, for all $t \leq t_0$, where $C$ is a
constant independent of $p$ and $t$.
\end{lemma}

\begin{proof}
$(B(p_i,t),g_i)$ converge to $(B(p,t),g_{\infty})$ in a complete
metric space $(M_{\infty},g_{\infty})$. Since the Bishop comparison
principle holds for orbifolds as well, we have that for every $\delta
\le t$:

\begin{equation}
\label{equation-comparison}
\frac{\vol_{g_i}B(p_i,t)}{t^4} \leq \frac{\vol_{g_i}B(p_i,\delta)}{\delta^4}
\end{equation}

$p_i$ is an orbifold point, with a curvature estimate 

\begin{equation}
\label{equation-equation_curvature_estimate}
|R(g_i)|(x) \leq \frac{\epsilon(r_i(x))}{r_i^2(x)}
\end{equation}

where $r_i(x) = \rho_i(x,p_i)$. Let $\Delta^*_r$ denote the punctured
ball in $C^2$ with radius $r$ and let $g_F$ be a standard euclidean
metric.

\begin{claim}
\label{claim-claim_help}
For any $i$ there exists $\delta_i > 0$ and a diffeomorphism $f_i$
from $\Delta^*_{\delta_i}$ into the universal covering $E_i$ of
$B(p_i,\delta_i)$ such that the covering map $\pi_i: E_i \to
B(p_i,\delta_i)$ is finite and

$$\max_{\Delta^*_{\delta_i}}|(\pi_i\circ f_i)^*g_i - g_F|_{g_F} \le
\epsilon_i$$

\end{claim}

\begin{proof}{(sketch)}

Call singular points $x_{\infty\beta}$ singular points of type I and
$p_{\infty i}$ singular points of type II. The total number of
singular points in each $(M_i,g_i)$ (after taking a subsequence if
necessary) is $N+S$. Denote this number by $K$.

The proof of the claim \ref{claim-claim_help} is just a modified proof
of lemma $3.6$ in \cite{tian1990}. For the convenience of a reader we
will just give a sketch of a proof here.

Let $E_{\kappa}(r) = \{x\in M_{\infty}\:\: |\:\:
\rho_{\infty\kappa}(x) < r\}$, where $\kappa$ is one of the indices
$\beta$ or $i$ (depending on the type of a singularity). Shortly, we
will say that $0\le \kappa \le K$. $\rho_{\infty\kappa}(x)$ is a
distance from a singular point in consideration to a point $x$. The
same argument as in \cite{tian1990} (lemma $3.4$) tells us that there
is a constant $L$ indpendent of $r$ such that the number of the
connected components in $E_{\kappa}(r)$ is less than $L$ for any $1
\le \beta \le K$.

We have that $|R(g_i)|(x) \le \frac{\epsilon(r_i(x))}{r_i^2(x)}$. By
taking a limit on $i$ and using the definition of convergence we get
that the same inequality holds for a limit metric $g_{\infty}$.

Fix some orbifold $(M_i,g_i)$. Consider one of its singular points
$p_i$. We will ignore subscripts for a moment (keepning in mind that
we are on some orbifold of our original sequence). We will show that
for any $\epsilon \in (0,1)$ there is a $r_{\epsilon} > 0$ such that
for any $r > r_{\epsilon}$, there is a diffeomorphism $\phi_r$ from an
annulus $\Delta(\frac{r}{2}, 2r)\subset C^2$ into $\pi^{-1}D(r,2)$
with its image containing $\pi^{-1}(D(r,2 - \epsilon))$ and

\begin{equation}
\label{equation-equation_estimate1}
\max\{||\phi_r^*\pi^* g - g_F||_{g_F}(x) | x\in
\Delta(\frac{r}{2},2r)\} \le \epsilon
\end{equation}

where $D(r,2) = \{x\in M\:\: |\:\: \frac{r}{2} \le \rho(x) \le r\}$
($\rho(x)$ is a distance from $x$ to a singular point $p$ in
consideration).

We prove it by contradiction. If this is not true, there is a sequence
$\{r(j)\}$ with $\lim_{j\to\infty} r(j)\to 0$ such that for any $r(j)$
no diffeomorphism with the above property exists. Since $p$ is an
orbifold point with a structure group $\Gamma < SU(2)$, where
$|\Gamma|$ is uniformly bounded for our sequence of orbifolds by
theorem \ref{theorem-order-singularity}, by our curvature estimate
\ref{equation-equation_curvature_estimate} $(D(2,r(i)),
\frac{1}{r(i)^2}g)$ converge to $\Delta(\frac{1}{2},2)/\Gamma$ in
$C^2\backslash\{0\}/\Gamma$. Since the estimate
\ref{equation-equation_estimate1} is invariant under scaling, by the
definition of convergence we immediately get a contradiction.

At the end we just glue all $\phi_r$ together to obtain the required
local diffeomorphism $f_i$ in the statement of our claim
\ref{claim-claim_help}.

\end{proof}

For each $i$ choose $\delta_i$ as in claim \ref{claim-claim_help}. We
can assume that $\delta_i\to 0$ as $i\to\infty$ (by decreasing
$\delta_i$ if necessary). By the claim $B(p_i,\delta_i)$ is covered by
a smooth manifold $E_i$, with a covering group $\Gamma_i$ ( a subset
of $SU(2)$, since all our singular points are rational) such that the
smooth manifold is diffeomorphic to a ball $\Delta_{\delta_i} \in C^2$
of radius $\delta_i$ via diffeomorphism $f_i$, where:

\begin{equation}
\label{equation-equation_close_to_euclidean}
|f_{i}^*\pi_i^*g_{\infty} - g_F|_{g_F} \leq \epsilon_i
\end{equation}

where $\epsilon_i$ tends to $0$ when $i \to \infty$, $g_F$ is a
standard euclidean metric and $\pi_i$ is just a covering map. Then:

\begin{equation}
\label{equation-help}
\frac{\vol_{g_i}B(p_i,\delta_i)}{\delta_i^4} =
\frac{\vol_{(\pi \circ f_i)^*g_F}\Delta_{\delta_i}}{|\Gamma_i|\delta_i^4}
\end{equation}

where $\Gamma_i$ is bounded by a constant that does not depend on $i$
by theorem \ref{theorem-order-singularity}.

By estimate \ref{equation-equation_close_to_euclidean} we have

$$\lim_{\delta_i\to 0}\frac{\vol_{(\pi \circ
f_i)^*g_F}\Delta_{\delta_i}}{\delta_i^4} = w_n$$

where $w_n$ is a volume of a unit euclidean ball. Letting $\delta_i\to
0$ in \ref{equation-comparison}, we get that:

$$\vol_{g_i}B(p_i,t) \leq Ct^4$$

We will get the result letting $k \to \infty$ and $i \to \infty$ in
the inequality above.

\end{proof}

So far we have proved that $\{x_{\infty\beta}\}_{1\le\beta\le N}$ are
the orbifold points of $M_{\infty}$ and points $\{p_{\infty
i}\}_{1\le\beta}$ have the following property: for any $p_{\infty i}$,
there is a neighbourhood $U_i$ of $p_{\infty i}$ in $M_{\infty}$ such
that any connected component $U_{ij}$ $(1 \le j \le l_i)$ of $U_i\cap
(M_{\infty}\backslash \sing(M_{\infty}))$ is covered by a smooth
manifold $\tilde{U}_{ij}$ with the covering group $\Gamma_{ij}$
isomorphic to a finite group in $U(2)$ and $\tilde{U}_{ij}$ is
diffeomorphic to a punctured ball $\Delta_{\tilde{r}}^*$ in
$\mathrm{C}^2$. If $\phi_{ij}$ is a diffeomorphism from
$\Delta_{\tilde{r}}^*$ onto $\tilde{U}_{ij}$ and $\pi_{ij}$ a covering
map from $\tilde{U}_{ij}$ onto $U_{ij}$, then the pull-back metric
$\phi_{ij}^*\circ\pi_{ij}^*(g_{\infty})$ extends to a smooth metric on
the ball $\Delta_{\tilde{r}}$, i.e. $g_{\infty}$ extends to a smooth
orbifold metric on each component of $U_{ij}\cap
M_{\infty}$. Therefore $(M_{\infty},g_{\infty})$ is a connected
Kahler-Einstein orbifold (maybe reducible) with finitely many singular
points.

To finish the proof of our main theorem in this section we have to
show that $(M_{\infty},g_{\infty})$ is locally irreducible, that is,
for any singular point $p_{\infty i}$ the punctured ball $B(p_{\infty
i},r)\backslash \{p_{\infty i}\}$ is connected for small values of
$r$.

We will call $M_{\infty}$ a generalized orbifold.

\begin{lemma}
\label{lemma-lemma_connect}
For some small $t_0$, $B(p,t)\backslash\{p\}$ is connected for all $t
\le t_0$.
\end{lemma}

We will postpone the proof of this lemma till the next section where
we will generalize some results about plurianticanonical sections to
the case of a sequence of Kahler-Einstein orbifolds and use these
results to prove lemma \ref{lemma-lemma_connect}.

\begin{lemma}
$(M_i,\tilde{g}_i,y_{i\beta})$ converges to
$(M_{\beta},h_{\beta},y_{\beta})$ in the pointed Gromov-Hausdorff
distance, where $(M_{\beta},h_{\beta})$ is a complete, non-compact,
Ricci flat,non-flat ALE orbifold with one end. $\tilde{g}_i =
\frac{1}{r(i)^2}g_i$, where $r(i)\to 0$ as $i\to\infty$ and
$\{y_{i\beta}\}$ for $0\le\beta\le K$ is a set of singular points of
$M_i$ and points where curvature operators concentrate.
\end{lemma}

\begin{proof}

By considering $A_k^{r_i} = \{x\in M_i\:\: |\:\: \frac{r_i}{k} \le
\rho_{\infty}(p_i,x) \le kr_i\}$ with metric $\tilde{g}_i$, standard
arguments as in \cite{anderson1989}, \cite{bando1989} and
\cite{tian1990} will give us that a sequence of pointed orbifolds
$\{(M_i, \tilde{g}_i, p_i)\}$ converge to a complete, non-compact,
Ricci-flat, non-flat $2$ compex orbifold which is ALE of order $3$
with 1 or more ends. Assume that it has 2 ends. Call it
$\tilde{M}$. The assumptions on our original sequence of orbifolds
give us a non-collapsing condition: $\vol_i B_i(x,r) \ge Cr^4$ for all
$i$ and all $x\in M_i$. It is invariant under scaling, so it will hold
also on our limit manifold $\tilde{M}$. Since $\tilde{M}$ has $2$
ends, it splits off a line and therefore $\tilde{M} = N\times R^k$ (by
a splitting theorem for orbifolds proved by J. Borzellino in
\cite{jborzell2000}). If $N$ were not a compact orbifold it would
contain 2 ends by assumption and therefore we could apply splitting
theorem to $N$ again. At the end we get that either $\tilde{M} =
N\times R^k$ where $N$ is a compact orbifold and $1 \le k \le 3$ or
$\tilde{M}$ is one of the following orbifolds: $R^4$, a product of
$R^3$ with a closed ray, or a product of $R^3$ with a closed
interval. In the former case, since $N$ is compact and therefore of a
finite volume we get a contradiction with a volume noncollapsing
condition for $\tilde{M}$. In the later case $\tilde{M}$ would be flat
which is not true.
  
\end{proof}
 
Now we can state the main theorem of this section which proof follows
immediatelly from what we have said and proved in a discussion above.

\begin{theorem}
$p$ is an orbifold point of the completion of $M_{\infty}$ that we
will call $\overline{M}_{\infty}$.
\end{theorem}

\section {Sections of plurianticanonical bundles of orbifolds}

In this section we want to generalize some results of \cite{tian1990}
about the sections of plurianticanonical bundles of a sequnce of
smooth surfaces to the sections of plurianticanonical bundles of a
sequence of $2$ dimensional orbifolds. We want to show that a sequence
of sections of plurianticanonical bundles of our orbifolds converge in
the sense that we will define below, to a section of a
plurianticanonical bundle of a limit orbifold. We also want to obtain
some estimates on the norms of the limits of the sections of the
plurianticanonical bundles of a limit orbifold.

Let's start with the following definition:

\begin{definition}

Let $S_i\in H^0(M_i,K_{M_i}^m)$ where $\{M_i,g_i\}$ is a sequence of
$2$ dimensional orbifolds as above. By results in the previous section
we may assume that $(M_i,g_i) \to (M_{\infty},g_{\infty})$. Let
$\phi_i$ be diffeomorphisms from the definition of convergence (
i.e. for any compact set $K\subset M_{\infty}\backslash
\sing(M_{\infty}$ there are diffeomorphisms $\phi_i$ from compact
subsets $K_i\subset M_i$ onto $K$ such that $(\phi_i^{-1})^*g_i$ and
$\phi_{i*}\circ J_i\circ(\phi_i^{-1})_*$ uniformly converge on $K$ to
$g_{\infty}$ and $J_{\infty}$, respectively. We will say that $S_i$
converge to $S_{\infty}$ if for any compact subset $K\in
M_{\infty}\backslash(\sing(M_{\infty})$ and $\phi_i$ as above, the
sections $\phi_{i*}(S_i)$ converge on $K$ to a section $S_{\infty}$ of
$K_{M_{\infty}}^{-1}$ in $C^{\infty}$ topology.

\end{definition}

Let $(M_{i},g_{i})$ be a sequence of two orbifolds as above. Then we
have the following lemma:

\begin{lemma}
\label{lemma-convergence}

Let $S^{i}\in H^{0}(M_{i},K_{M_{i}}^{-m})$ with
$\int_{M_{i}}||S^i||_{g_i}dV_{g_i} =1$.  Then $\exists$ a subsequence
$S^{i_k}$ that converges to $S^{\infty}\in H^0 (M_{\infty},
K^{-m}_{M_{\infty}})$.

\end{lemma}

\begin{proof}

If $f_i = ||S^i||^2_{g_i}$, then we have that:
$$\Delta_i f_i = ||D_i S^i||^2 - 2mf_i \geq -2mf_i$$

on $M_i\backslash \sing{M_i}$. Omit  subscript $i$ in a further
discussion. Then:

\begin{equation}
\label{equation_equation_main}
-\Delta f \leq 2mf 
\end{equation}

on $M\backslash \sing M$, where $p$ is an orbifold point of $M$.

To conlude that $Df\in L^2$ we have to use a fact that $f\in L^q$ for
all $q$. The proof of this fact can be found in \cite{bando1989}. We
can take a cut-off function $\phi$ so that $0\le \phi \le 1$, $\phi =
0$ in $B(r')\cup (M\backslash B(2r))$, $\phi = 1$ on $B(r)\backslash
B(2r')$ with $|D\phi|\le Cr'^{-1}$ for $2r'\le r$. Then since $f\in
L^q$ for any $q$, after multiplying the inequality
\ref{equation_equation_main} by $\phi f$ and performing a partial
integration, we get:

$$\int_M \nabla f\nabla(f\phi) \le C\int_M f^2\phi$$ 

\begin{eqnarray}
\label{equation-equation_gradient_estimate}
\int_M |\nabla f|^2\phi &\le& C + \int_M |\nabla\phi|f|\nabla f|
\nonumber \\ &\le& C + \tilde{C}\int_M|\nabla\phi|^2 f^2 +
\epsilon\int_M |\nabla f|^2 \nonumber \\ &\le& C + \tilde{C}(\int_M
|\nabla \phi|^{2p})^{p^{-1}} (\int_M f^{2q})^{q^{-1}} + \epsilon\int_M
|\nabla f|^2
\end{eqnarray}

where we have used Cauchy-Schwartz and Holder inequalities with
$\epsilon < 1$ and $p < n$. We let $r'\to 0$ and $r\to \diam{M}$ in
\ref{equation-equation_gradient_estimate} to get that $Df\in L^2$.

Assume without a loss of generality that $p$ is the only singular point
of $M$.

\begin{claim}
$\forall \eta \in C_0^{0,1}(B(p,r))$, where p is a singular point of
$M$, $r > 0$ arbitrary  and $B=B(p,r)\backslash \{p\}$: 

$$\int_{B} D\eta Df \leq 2m\int_B \eta f$$.
\end{claim}

\begin{proof}

Let $\eta_k \in C^1_0(B\backslash \{p\})$ s.t. $\eta_k \to 1$ a.e. and
$\int_B|D\eta_k|^2 \to 0 \ \ (k\to \infty)$.

\begin{eqnarray*}
\int D\eta \eta_k Df &=& \int D(\eta \eta_k) Df - \int D\eta_k \eta Df \nonumber \\
&\le&  2m\int \eta\eta_k f + (\int |D\eta_k|^2)^{\frac{1}{2}}(\int
\eta^2|Df|^2)^{\frac{1}{2}}
\end{eqnarray*}

Let $k \to \infty$. Then

$$\int_B D\eta Df \leq 2m\int_B \eta f \ \ \ \ \ \forall \eta \in C_0
^{0,1}(B(p,r))$$

\end{proof}

\begin{claim}

$\forall \eta \in C^1(M)$,

$$\int_{M} D\eta Df \leq 2m\int_{M} f\eta$$.

\end{claim}

\begin{proof}

Take $U_k = B(p,r + \frac{1}{k})$ and $V_k=M\backslash B(p,r -
\frac{1}{k})$ to be the open covering of $M$ and let $\phi_k$,
$\psi_k$ be the partititon of unity subordinated to $U_k,V_k$.  By
using the previous claim and the fact that $V_k$ is smooth we have:

\begin{eqnarray*}
\int_M D\eta dF &=& \int_M D(\eta \phi_k +\eta \psi_k)Df  \\
&=& \int_{U_k}D(\eta\phi_k)Df + \int_{V_k}D(\eta \psi_k)Df\\ 
&\le&  (2m\int_{U_k}f\eta\phi_k + 2m\int_{V_k}f\eta\psi_k) \\
&\le& 2m\int_M f\eta(\phi_k + \psi_k) 
= 2m\int_{M}f\eta
\end{eqnarray*}

\end{proof}
 
Now by Moser's iteration argument and Sobolev lemma we get that:

$$\sup_{M_i \backslash \sing (M_i)}||S^i||_{g_i}(x) \leq C(m) \ \ \ \ \
\forall i$$ 

Similarly like in the case of a sequence of smooth KE surfaces with
positive first Chern class, to finish the proof of lemma
\ref{lemma-convergence} one can prove, using the Caushy integral
formula that the $l$th covariant derivatives of $\Phi_{i*}S^i$ are
uniformly bounded on compact sets $K_i \in M_i\backslash \sing M_i$ ,
by a constant depending on $l$ and $K=\Phi_i(K_i)$, where $\Phi_i$ are
diffeomorphisms from the definition of convergence of a sequence of
orbifolds. Since $(\Phi^{-1})^*g_i$ uniformly converge to $g_{\infty}$
in $K$, the lemma is proved. 

\end{proof}

Analogously like in \cite{tian1990} for the smooth case, it can be
shown that if $S\in H^0(M_{\infty}, K_{M_{\infty}}^{-m})$, $\exists$ a
sequence $S^i\in H^0(M_i, K_{M_i}^{-m})$ converging to $S$.  If we
prove that $h^0(M_j,K_{M_j}^{-m})$ is bounded above uniformly in $j$,
by taking a subsequence (denote it again by $\{M_j\}$),
$h^0(M_j,K_{M_j}^{-m})=h^0(M_{\infty},K_{M_{\infty}}^{-m})$, which
implies that  $\{S_{\beta}^i\}_{0\leq\beta\leq N_m}$, an
orthonormal basis of $H^0(M_i, K_{M_i}^{-m})$  converges to an
orthonormal basis of $H^0(M_{\infty}, K_{M_{\infty}}^{-m})$.

\begin{lemma}
\label{lemma-lemma_boundness_of_PAD}
$h^0(M_i,K_{M_i}^{-m})\leq C(m),\ \ \ \forall i$ , where $C(m)$ is a
constant that depends only on $m$.
\end{lemma}

\begin{proof}
Omit subscript $i$ in the proof of the theorem. Let $M$ be an orbifold
with $N$ singular points $p_1,\dots,p_N$ of orders
$|\gamma_1|,\dots,|\gamma_N|$, bounded uniformly by $C$. By the
generalized Rieman-Roch formula for an orbifold $M$, we have that:

\begin{equation}
\label{equation-equation_generalized}
\chi(M,K_M^{-m}) = \frac{1}{2} \int_M m^2c_1(M)^2 + \frac{1}{2}\int_M
mc_1(M)\pi_*c_1(\tilde{M}) + \chi(O_M)
\end{equation}

where $\pi:\tilde{M}\to M$ is a resolution of singularities.  Because
of the uniform bounds on $c_1^2(M_i)$ and $c_2(M_i)$ of our sequence
of orbifolds $(M_i,g_i)$ at the beginning, after applying formula
\ref{equation-equation_generalized} to $(M_i,g_i)$ we get that:

$$(h^0 - h^1 + h^2)(M, K_M^{-m}) \leq C(m)$$ 
where $C(m)$ is a constant independent of $i$.

Since $K_M^{-m} \ge 0$ by assumption, it follows that $h^1 = h^2 =0$
by Bailey's version of Kodaira's vanishing theorems for
orbifolds. Now the lemma follows. 

\end{proof}

We also have the following theorem:

\begin{theorem}
\label{theorem-comparison}
Let $(M_i,g_i)$ be a sequence of $2$ dimensional Kahler-Einstein
orbifolds as above, converging to a Kahler-Einstein orbifold
$(M_{\infty},g_{\infty})$ such that $h^0(M_i,K_{M_i}^{-m}) =
h^0(M_{\infty}, K_{M_{\infty}}^{-m})$ ( we can assume this by the
previous lemma). Let $\{S_{m\beta}^i\}_{0\leq\beta\leq N_m}$ be a
sequence of linearly independent sections of $K_{M_i}^{-m}$ and let
$\{S_{m\beta}^{\infty}\}_{0\leq \beta \leq N_m}$ be a basis of
$H^0(M_{\infty}, K_{M_\infty}^{-m})$. Then:

$$\lim_{i\to \infty}(\inf_{M_i}\{\sum_{\beta
=0}^{N_m}\||S_{m\beta^i||^2 }) \geq \inf_{M_{\infty}}\sum_{\beta
=0}^{N_m}||S_{m\beta}^{\infty}||^2_{g_{\infty}}$$.

\end {theorem}

\begin{proof}

$$\Delta_i||D_iS_{m\beta}^i||^2_{g_i} =
||D_iD_iS_{m\beta}^i||^2_{g_i} - (4m - 1)||D_iS_{m\beta}^i||^2_{g_i}
\geq  -(4m - 1)||D_iS_{m\beta}^i||^2_{g_i}$$

Like in lemma \ref{lemma-convergence} we can get that:

$$\sup\{||D_iS_{m\beta}^i||^2_{g_i}(x)\:\: |\:\: 0\leq \beta \leq
N_m,\:\: x\in M_i\backslash \sing(M_i)\} \leq C_1(m)$$

Combining this with the result of lemma \ref{lemma-convergence}, we
get that the first derivatives of $f_i =
\sum_{\beta=0}^{N_m}||S_{m\beta}^i||^2_{g_i}$ are uniformly bounded on
$M_i\backslash \sing(M_i)$. Let $L_i = \inf_{M_i}f_i$ and $\epsilon >
0$. Then $\exists x_i \in M_i$, s.t. $L_i > f_i(x_i) -
\epsilon$. Since $|Df_i|\leq C$, it follows that $(*)\ \
\omega_i(f_i,r)\leq Cr,\ \ \ \forall i$, where $w_i(f_i, \cdot)$ is an
oscilation of $f_i$. Take $\{r_i\}$ s.t. $r_i \to 0 \ \ (i\to
\infty)$and $Cr_i < \epsilon \ \ \forall i$. Let $z_i \in \partial
B_{\frac{r_i}{2}}(p_i^j)$, where $p_i^j$ is a singular point in
$M_i$. We will have 2 cases:

\begin{enumerate}

\item

for almost all $i$ $x_i \in B(p_i^j,r_i)$:

$$L_i \ge f_i(z_i) - 2\epsilon = (\phi_i^{-1})^*f_i(\phi_i(z_i)$$

because of $(*)$ when $r = r_i$, where $\phi_i: M_i\backslash
\bigcup_{j=1}^N B_{\frac{r_i}{4}}(p_i^j,g_i)\to M_{\infty}$ are the
embeddings as in a definition of convergence and $N$ is the number of
singular points on $M_i$ (by taking a subsequence, we may assume that
$N$ does not depend on $i$).

Let $i\to \infty$. Since $Im \phi_i$ tends to $M_{\infty}\backslash
\sing M_{\infty}$ and $f_i \to f$ by lemma \ref{lemma-convergence}, we
have that:

$$\lim_{i\to \infty} L_i \geq \inf_{M_{\infty}} f - 2\epsilon$$

Since the above inequality holds for every $\epsilon > 0$, the
statement is proved.

\item

$x_i \notin B(p_i^j,r)$ for infinitely many $i$, for some $r > 0$
fixed (if such $r$ did not exist, we would find a sequence $\{r_i\}$
with the properties like in the case $(1)$). The statement follows by
the same arguments as for the case of a sequence of smooth surfaces 

(look at \cite{tian1990}).

\end{enumerate}

\end{proof}

Let $\mathcal{A}(C_1,C_2)$ denote the same set of orbifolds as at the
beginning. Let $\mathcal{A}_n$ be the set of all orbifolds $M\in
\mathcal{A}(C_1,C_2)$, such that $h^0(M,K_M^{-1}) = n$, where $n \le
C(1)$, where $C(1)$ is a constant as in lemma
\ref{lemma-lemma_boundness_of_PAD}.

Since we have theorem \ref{theorem-comparison}, the same arguments as
in \cite{tian1990} give us the following estimate:

\begin{theorem}

There are a universal integer $m_0 > 0$ and a universal constant $C >
0$, s.t. for any KE orbifold $(M,g)$ in $\mathcal{A}_n$  we have:

$$\inf_M \{\sum_{\beta=0}^{N_{m_0}} ||S_{\beta}||^2_g\} \ge C > 0$$
where $N = N_m+1$ is the complex dimension of $H^0(M,K_M^{-m_0})$, and
$\{S_{\beta}\}_{0\le \beta \le N}$ is an orthonormal basis of
$H^0(M,K_M^{-m_0})$, with respect to the inner product induced by $g$.
\end{theorem}

\begin{proposition}
The generalized Kahler-Einstein orbifold $(M_{\infty},g_{\infty})$
that we constructed in the previous section is locally irreducible,
i.e. for every $p\in M_{\infty}$ there exists some $r>0$ such that
$B(p,r)\backslash\{p\}$ is connecetd.
\end{proposition}

\begin{proof}
For the proof we will refer to \cite{tian1990} (see proposition
5.2. in \cite{tian1990}).
\end{proof}

\end{document}